\documentclass[11pt,a4paper]{article}
\usepackage{latexsym}

\title{Graphs, flags and partitions}

\author{Jonathan Fine\relax
\thanks{203 Coldhams Lane, Cambridge, CB1 3HY, England.  
\quad E-mail: \texttt{j.fine@pmms.cam.ac.uk}}
}
\date{17 September 1998}

\newtheorem{theorem}{Theorem}
\newtheorem{definition}[theorem]{Definition}
\newtheorem{proposition}[theorem]{Proposition}
\newtheorem{lemma}[theorem]{Lemma}
\newenvironment{proof}
  {
    \addvspace {3pt}%
    \noindent\emph{Proof.~~}%
  }
  {
    \nobreak\hfill$\Box${\parfillskip 0pt \par}%
    \addvspace {6pt}%
  }

\newcommand\lambdabar{\overline{\lambda}}
\newcommand\Gbar{\overline{G}}
\newcommand\Gcal{{\cal G}}
\newcommand\Zcal{{\cal Z}}
\newcommand\fbar{\overline{f}}  %%  use smash and phantom?

\newcommand\bibrule{\rule{2pc}{0.4pt}}

\begin{document}
\maketitle
\begin{abstract}
\noindent
This paper defines, for each graph $G$, a flag vector $fG$.  The flag
vectors of the graphs on $n$ vertices span a space whose dimension is
$p(n)$, the number of partitions on $n$.  The analogy with convex
polytopes indicates that the linear inequalities satisfied by $fG$ may be
both interesting and accessible.  Such would provide inequalities both
sharp and subtle on the combinatorial structure of $G$.  These may be
related to Ramsey theory.
\end{abstract}

\section{Introduction}

The formulation and proof of subtle combinatorial inequalities on the
number of vertices, edges and so forth on a graph is the goal towards
which this paper is intended to be the first step.  This paper puts
forward a precise definition of what is meant by `and so forth'.  More
exactly it assigns to each graph $G$ a flag vector $fG$.  The span of
$fG$, as $G$ ranges over all graphs on $n$ vertices, has dimension equal
to $p(n)$, the number of partitions on $n$.  The linear inequalities
satisfied by $fG$, or what comes to the same thing the convex hull
$\Delta(n)$ of $fG$ as $G$ ranges over all $n$-vertex graphs, seems not
to be arbitrary, at least for small values of $n$.

Ramsey theory can be expressed in terms of combinatorial inequalities. 
For example, that any graph $G$ on six vertices contains a $3$-element
subset $S$ of its vertex set, such that the restriction $G|_S$ of $G$ to
$S$ is either the complete graph or the empty graph:  this statement is
equivalent to the following.  For each $3$-element subset $S$ of the
vertex set, temporarily let $\lambda(G,S)$ be $1$ if $G|_S$ is the
complete graph, and zero otherwise.  Now let $\lambda(G)$ be the sum of
$\lambda(G,S)$ over all $S$.  Define $\lambdabar(G)$ to be as
$\lambda(G)$, but this time counting empty graphs.  The Ramsey statement
is equivalent to $\lambda(G) + \lambdabar(G) \geq 1$.  Although this
equivalence is trivial, it does indicate the sort of relationship we are
looking for.  Whether or not Ramsey-type statements can be formulated in
terms of the flag vector $fG$ is not yet clear.  If so, they may very well
involve non-linear functions of $fG$.

The flag vector can be defined in several different ways, that are
essentially equivalent in the sense that each is a linear function of any
of the others.  The easiest, although lacking in motivation, is via
subgraphs.  Here is its definition.

\begin{definition}
Suppose $G$ is a graph on $n$ vertices, and that $H$ is a subgraph of
$G$, also on $n$ vertices.  Each such subgraph will contribute
$\lambda^s(H)$, a quantity defined below.  In other words
\[
    f^sG = \sum \nolimits _ { H } \lambda ^s ( H ) 
    \qquad \mbox { $H$ a subgraph of $G$}
\]
defines the \emph{subgraph form $f^sG$} of the flag vector.
\end{definition}

\begin{definition}
The \emph{subgraph contribution $\lambda^s(G)$} is defined as follows. 
Each graph $G$, via connectivity, determines a partition of its vertex
set, which can then be thought of as a partition $\pi=\pi(G)$ of the order
$n$ of the graph.  The quantity $\lambda^s(G)$ will be a formal sum of
partitions, and the contribution due to $G$ will be $s^a(G)\pi(G)$, where
$s^a(G)$ is next to be defined.
\end{definition}

\begin{definition}
Let $G$ be a graph on $n$ vertices.  A \emph{shelling} $\sigma$ of $G$ is
a sequence $G=G_0$, $G_1$, $\ldots$, $G_n$ of graphs obtained from $G$ by
removing the vertices from $G$ one at a time.  As a vertex is removed, so
all edges through that vertex are also removed.  An \emph{acyclic
shelling} is a shelling where at most one edge is removed at each stage. 
The \emph{acyclic shelling number} $s^a(G)$ of $G$ is the number of
acyclic shellings that $G$ possesses.  For example, if $G$ has a cycle then
$s^a(G)$ is zero.
\end{definition}

The flag vector is an example of a rule $v$ that assigns to each graph $G$
a vector $v(G)$ lying in some vector space $V$.  For such to be useful it
should retain some, but perhaps not all, of the geometric information in
$G$.  For example, say that $v$ \emph{distinguishes graphs} if
$v(G)=v(G')$ implies that $G$ and $G'$ are isomorphic as graphs.  The
location of $G$ in an enumeration of all isomorphism classes of graphs is
an example of such a rule.  The essentially random ordering ordering of
graphs that such provides is unlikely to be useful.  Similarly, say that
$v$ \emph{embeds graphs} if the vectors $v(G)$, as $G$ ranges over all
isomorphism types of graphs, are linearly independent.  This is even
worse.  Enumeration gives only random relations between the vectors
$v(G)$.  Embedding gives none at all.

Much of the motivation for this paper has come from the theory of convex
polytopes.  In this case the definition of the flag vector is
straightforward, and many subtle combinatorial inequalities are either
known or conjectured.  When $\Delta$ is a simple $n$-polytope (exactly $n$
facets at each vertex) McMullen's conjectured conditions on the face
vector \cite{bib.McM.NFSP} were proved sufficient (by Stanley
\cite{bib.RS.NFSP}) and necessary (by Billera and Lee
\cite{bib.LB-CL.SMC}). In the general case the `mpih $h$-vector' is known
to be unimodal when $\Delta$ has rational vertices
\cite{bib.JD-FL.IHNP,bib.KF.IHTV,bib.RS.GHV}.  Latent in the author's
local-global variant of intersection homology
\cite{bib.JF.LGIH,bib.JF.CPLA,bib.JF.SSIH} and his analogue of ring
structure \cite{bib.JF.RSUE} are various conjectures on polytope flag
vectors.

To provide some background, we here recall the basic definition and main
result of \cite{bib.MB-LB.gDS}.  (For the subtle inequalities, see the
works just cited.)

\begin{definition}
Suppose $\Delta$ is a simple $n$-dimensional polytope.  In that case the
\emph{face vector} $f\Delta$ of $\Delta$ is the sequence $f\Delta=(f_0,
\dots, f_{n-1})$, where $f_i$ is the number of $i$-dimensional faces on
$\Delta$.  Now suppose that $\Delta$ is a general $n$-dimensional
polytope.  A \emph{flag} $\delta$ of $\Delta$ is a sequence
\[
    \delta =
      ( \delta_1 \subset \delta_2 \dots \subset \delta_r \subset \Delta )
\]
of faces, each strictly contained in the next.  Its \emph{dimension
sequence}, or \emph{dimension} for short, is the sequence
\[
    d = ( d_1 < d_2 < \dots < d_r < n )
\]
of integers, where $d_i$ is the dimension of $\delta_i$.  The \emph{$d$-th
component} $f_d\Delta$ of the \emph{flag vector} $f\Delta$ is the number
of $d$-dimensional flags on $\Delta$.  One can think of $d$ as a subset
$I$ of $\{0,1,\dots,n-1\}$, and so $f\Delta$ has $2^n$ components, indexed
by such subsets.  (If $\Delta$ is known to be simple, then there is a
linear function that computes the flag vector from the face vector.  Thus,
the overloaded use of $f\Delta$ will not in this situation do any harm.)
\end{definition}

\begin{theorem}[Generalised Dehn-Sommerville equations]
The span of the flag vectors of $n$-dimensional polytopes has dimension
the $(n+1)$st Fibonacci number
\end{theorem}

\section{The verbose flag vector $f^vG$}

This definition takes its motivation from convex polytopes.  Suppose
$\Delta$ is an $n$-dimensional convex polytope.  To define the flag vector
of $\Delta$ it is enough to know what it means, for $\delta$ to be a
facet of $\Delta$.  For once this is known, the faces are facets of facets
of facets etc., and the flag vector follows from the counting of chains.

What then is a `facet' of a graph $G$ on $n$ vertices?  Clearly, it should
be a graph $H$ on $n-1$ vertices.  Such can be obtained by removing a
single vertex $v$ from $G$, and when this is done some edges may have to
be removed as well.  The flag vector will record such information.

\begin{definition}
Suppose $v$ is a vertex of a graph $G$.  As usual, the \emph{order} $m_v$
of $v$ (on $G$) is the number of edges of $G$ that pass through $v$. 
Define the $v$-facet $G_v$ of $G$ to be the result of removing $v$ (and
also the $m_v$ edges) from $G$.
\end{definition}

\begin{definition}
The \emph{verbose flag vector} of the unique graph on zero vertices is the
number $1$.  For other graphs the \emph{verbose flag vector} $f^vG$ is a
formal sum of words in $a$ and $b$ obtained by the following recursion:
\[
    f^vG = \sum \> (a+m_vb) f^v G_v
\]
where the sum is over all vertices $v$ (or facets $G_v$) of $G$.  If $G$ is
of order $n$ then $f^vG$ is a homogeneous polynomial of degree $n$ in the
non-commuting variables $a$ and $b$.  There are $2^n$ such words $w$ in
$a$ and $b$.  Each can be thought of as a subset $I$ of $\{1,\dots,n\}$.
\end{definition}

By way of example, let $3_0$, $3_1$, $3_2$ and $3_3$ denote the graphs on
$3$ vertices with $0$, $1$, $2$ and $3$ edges respectively.  The reader is
asked to verify the following results.  Later, we will give
Proposition~\ref{order.3.eqn} a geometric meaning.

\begin{proposition}
\label{order.3.f}
The verbose flag vectors of the graphs on $3$ vertices are as follows:
\[
\begin{array}{l}
    f^v(3_0) = 6 aaa \\
    f^v(3_1) = 6 aaa + 2 aba + 4 baa \\
    f^v(3_2) = 6 aaa + 4 aba + 8 baa + 4 bba \\
    f^v(3_3) = 6 aaa + 6 aba + 12 baa + 12 bba
\end{array}
\]
\end{proposition}

\begin{proposition}
\label{order.3.eqn}
There is a single linear relation among these flag vectors, namely
\[
    f( 3_0 ) - 3 f (3_1 ) + 3 f(3_2) - f(3_3) = 0 \>.
\]
\end{proposition}

Each form of the flag vector has its own advantages.  Here is a result
that is most easily proved by using the verbose form.  But first a
definition.

\begin{definition}
As usual, if $G$ is a graph on $n$ vertices, let the \emph{complementary
graph} $\Gbar$ be the graph on the same $n$ vertices, whose edges are
exactly the non-edges of $G$.  Let $\fbar G$ denote the flag vector (in
whatever form) of the complement $\Gbar$ of $G$.
\end{definition}

\begin{theorem}
The verbose flag vector of $\fbar^vG$ of the complement $\Gbar$ of $G$ is
a linear function of $f^vG$.  (This statement means that there is a linear
function $L$ on degree $n$ words in $a$ and $b$, such that the equation
$L(f^vG)=\fbar^vG$ holds true, for every $n$-vertex graph $G$.)
\end{theorem}

\begin{proof}
If $G$ has no vertices the result is trivial.  We will use induction.  The
$v$-contribution to $\fbar G$ is $(a + (n-1-m_v)b) \fbar G_v$, where $m_v$
is the multiplicity of $v$ on $G$.  Thus, the transformation
\[
    a \mapsto a + (n-1) b \qquad b \mapsto -b
\]
applied to the first letter of each word in $fG$, together with the
inductive formula applied to the remainder of the word, will define the
order $n$ rule that takes $f^vG$ to $\fbar^vG$.
\end{proof}

In fact the induction in this proof is spurious.  The transformation
\[
    a \mapsto a + (i-1) b \qquad b \mapsto -b
\]
applied to each of the letters $w_i$ of a word $w=w_nw_{n-1}\dots w_1$
will transform $fG$ to $\fbar G$.  It is easily seen to be an involution,
as it should be.

Shellings can be used to compute the verbose flag vector.  The next
theorem follows immediately.

\begin{definition}
Let $\sigma$ be a shelling of a graph $G$.  Define the \emph{verbose
shelling contribution} $\sigma^vG$ to $f^vG$ to be the product $w_1 w_2
\dots w_n$, where $w_i$ is $(a+m_ib)$, and where $m_i$ is the number of
edges of $G$ of the form $v_iv_j$ for $i<j$, and where $v_i$ and $v_j$ are
the $i$-th and $j$-th vertices in the shelling $\sigma$.
\end{definition}

\begin{theorem}
The formula $f^vG= \sum _ \sigma \sigma^v G$ computes the verbose flag
vector.
\end{theorem}

This gives another way of looking at $f^vG$.  First fix a word $w$ in $a$
and $b$.  Each shelling $\sigma$ will contribute to the $f_w$ component of
$f^vG$ as follows.  Whenever $w_i$ in $w=w_1\dots w_n$ is $b$, compute
$m_i$ as above.  The product of these $m_i$ is the $\sigma$ contribution
to the coefficient of $w$ in $f^vG$, or in other words to the $f_w$
component.

In Ramsey theory, non-constructive proofs of the following type are common.
First produce an overestimate of the number of graphs that do not have
some property $P$.  Next, count the number of graphs.  If, for some number
$n$ of vertices the first number is less than the second, then we have a
non-constructive proof that there is a graph that has the property $P$. 
This is an example of the `random method' that pervades graph theory and
combinatorics.  It is similar to the application of the pigeon-hole
principle.  

Here the property $P$ is a function $P(G)$ that takes the values zero and
one on graphs.  Now suppose that $l(G)$ is a non-negative linear function
of the flag vector $fG$.  (To produce such is to know something of the
linear inequalities on the flag vector.)  If we can compute the average
value of $l$ over all graphs (on $n$ vertices), and this value is
non-zero, then similarly we have a non-constructive proof that there is a
$G$ for which $l(G)$ is non-zero.

The earlier presentation of the verbose flag vector as a sum over all
shellings of a graph leads to the next result, whose proof is left to the
reader.

\begin{definition}
The \emph{total flag vector} $f(n)$ of all graphs on $n$ vertices is the
sum of $fG$, as $G$ ranges over all graphs on $n$ vertices.  Here, such a
graph is simply a subset of the $N={n \choose 2}$ possible edges on an
$n$-element vertex set.  Thus, there are $2^N$ such graphs.
\end{definition}

\begin{theorem}
For each word $w=w_n \dots w_1$, the $w$ component $f_w(n)$ is the product
\[
    n! \, \lambda_n \dots \lambda_1
\]
where $\lambda_i$ if $2^{i-1}$ is $w_i$ is $a$, while $\lambda_i$ is
$(i-1)2^{i-2}$ if $w_i$ is $b$.  In particular, if $l$ is a linear
function of $f$, then the average value of $l$ over all graphs can readily
be computed.
\end{theorem}

Given a word $w$, certain pairs of shellings $\sigma$ and $\sigma'$ will
make equal contributions to $f_w$, for any choice of an edge set for $G$.
Here is an example.  Let $w$ be $aabbaaa$, and suppose that $\sigma$ and
$\sigma'$ differ, if at all, only in the order in which the first two
vertices are presented, and in the order of the last three.  In other
words, we may permute $\{v_1,v_2\}$ and also $\{v_5,v_6,v_7\}$ to obtain
$\sigma'$ from $\sigma$.  Clearly there are 12 such possibilities and as
such changes do not affect $m_3$ and $m_4$, all these $\sigma'$ make the
same contribution to $f_w$.

From this it follows that for $w=aabbaaa$, the $f_w$ component of $f^vG$
is always divisible by $12$.  (Proposition~\ref{order.3.f} provides
similar examples.) This is a redundancy of sorts, which can be reduced by
using `partial shellings' that `respect' the word $w$.  There is another
redundancy in $f^vG$, which is that although it lies in a vector space of
dimension $2^n$, it does not span this space.  The next section will
define a flag vector that provides an integral spanning set for the space
in which it lies.

\section{The concise flag vector $f^cG$}

This section introduces the concise flag vector $f^cG$, which records the
same information as the verbose flag vector $f^vG$, but concisely.
In particular, like the subgraph flag vector $f^sG$, it takes
values in a space whose dimension is $p(n)$, the number of partitions of
the number of vertices on $G$, rather than the $2^n$ dimensional space
used by $f^vG$.

First we introduce a device that help us to understand
Proposition~\ref{order.3.eqn}.  It allows us to produce formal sums of
graphs, whose total flag vector is zero.  If we think of a graph as a
bunch of points joined by lines, then a \emph{graph with optional edges}
is the same, but some of the edges are represented by dashed lines. 
Between any two vertices there is at most one edge, either dashed or
ordinary.  Each dashed edge can be either completed or removed, and so
from a graph with $r$ dashed edges a total of $2^r$ graphs can be
produced.  The \emph{flag vector} $fG$ of a graph $G$ with optional edges
is the alternating sum of the flag vectors of the $2^r$ graphs so
obtained.  A graph with no optional edges is, of course, just an ordinary
graph.  The following formal definition takes care of the details.

\begin{definition}
A \emph{graph $G$ with optional edges} consists of a vertex set $V$,
together with disjoint sets $E$ and $C$ of unordered pairs of vertices.
Here $E$ is the regular edge set, while $C$ (for choice) is the set of
optional edges.  Such is usually regarded as a formal sum
\[
    (V,E,C) = 
      \sum \nolimits _ { B \subseteq C} \> (-1) ^ {|C|-|B|} \> (V, E\cup B)
\]
of normal graphs.  The quantity
\[
    f(G) =
      \sum \nolimits _ { B \subseteq C} \> (-1) ^ {|C|-|B|} \> f(V, E\cup B)
\]
is its \emph{flag vector}, of whatever type.
\end{definition}

Each graph with only optional edges is of course an alternating sum of its
subgraphs, thought of as ordinary graphs.  Conversely, each ordinary graph
can be written as a formal sum of graphs with only optional edges.  We can
think of an optional edge as the difference between an edge and a
non-edge.  Conversely, an edge can be thought of as the sum of an optional
edge and a non-edge.  The proof of the next result is left to the reader. 
It is at root nothing more than $x=(x-y)+y$.

\begin{theorem}
Let $G$ be a graph $(V,E)$.  For each subset $C$ of $E$, let
$(V,\emptyset,C)$ be the graph with $C$ as its optional edge set, and no
actual edges.  Then
\[
    (V,E) = \sum \nolimits _ { C \subseteq E } \> ( V, \emptyset , C )
\]
or in other words, $G$ is equivalent to a formal sum of graphs with only
optional edges, namely those associated to all its subgraphs.
\end{theorem}

The advantage of graphs with only optional edges is that their flag
vectors are much easier to calculate and otherwise deal with.  Applied to
the graph on $3$ vertices, with $3$ optional edges, the next result
establishes Proposition~\ref{order.3.eqn}.

\begin{theorem}
Suppose $G$ is a graph with optional edges, and in those edges a
cycle can be found.  Then the verbose flag vector $f^vG$ of $G$ is zero.
\end{theorem}

\begin{proof}
Let $w$ be a word in $a$ and $b$, and let $\sigma$ be a shelling of $G$. 
Let $v_i$ be the first vertex removed by $\sigma$ that lies on the
optional cycle. At this point $w_i$ is either $a$ or $b$.  By assumption,
at $v_i$ there are at least two optional edges of the form $v_i v_j$, for
$j>i$ in the $\sigma$ ordering.

If $w_i$ is $a$ and there is just a single optional edge $v_i v_j$, then
the $w$-$\sigma$ contribution to $f^vG$ is zero.  This is because changing
the choice made for $v_i v_j$ changes the sign of the contributions made
at this point.  If there are several optional edges, choose just one of
them and apply this argument to it.

If $w_i$ is a $b$, we will need two optional edges, say $v_i v_j$ and
$v_i v_k$.  This gives four $(2\times 2)$ choices altogther, and when we
determine their contributions we find that they are $+2$, $-1$ and $-1$,
and $0$.  (Both optional edges present, just one optional edge, and no
optional edge.)  The sum is zero and so the result follows.  As before, if
there are more than two optional edges, the result still holds, by
applying the argument to just two of them.
\end{proof}

The argument we just used can also be applied to compute the verbose flag
vector $f^vG$ of a graph $G$ with only optional edges.  To begin with, we
will apply it just to a tree.

\begin{theorem}
Suppose $T$ is a tree (connected acyclic graph) on $n$ vertices.  Consider
$T$ as a graph with only optional edges.  Then $f^v(T)$ is
$s^a(T)b^{n-1}a$, where $s^a(T)$ is the acyclic shelling number of $T$.
\end{theorem}

\begin{proof}
Let $\sigma$ be a shelling of $T$.  If $\sigma$ is not acyclic then the
previous argument, applied to the first vertex $v_i$ at which $\sigma$
fails to be acyclic, shows that the $\sigma$ contribution to $f^vT$ is
zero.  This holds whether $w_i$ in $w=w_1\dots w_n$ is $a$ or $b$.  Thus,
only acyclic shellings contribute.  Now suppose $w$ contains an $a$, and
that $w_i$ is the first $a$ in $w$.  Assume that $i<n$.  Provided $\sigma$
is acyclic, when $w_i$ is removed an optional edge will be removed also,
and the previous argument shows that this $w$-$\sigma$ contribution is
zero.  Finally, if $w_n$ is $b$ then the $w$ contribution must be zero,
for at the last vertex there is no edge left to remove.  Thus, the only
non-zero contributions are when $\sigma$ is an acyclic shelling, and then
$b^{n-1}a$ is the contribution.  By definition, $T$ has $s^a(T)$ such
shellings.
\end{proof}

Now suppose that $G = T_1 \sqcup T_2 \dots \sqcup T_r$ is an acyclic
graph, broken up as a disjoint union of trees.  As before, we will treat
$G$ as a graph with only optional edges.  Suppose that an acyclic
shelling $\sigma_i$ is given for each component $T_i$ of $G$.  From this
information many acyclic shellings $\sigma$ of $G$ can be constructed,
such that when $\sigma$ is restricted to each $T_i$ (and also renumbered
to account for absent vertices), the result is $\sigma_i$.  Looked at
abstractly, we have $r$ disjoint sets with an ordering $\sigma_i$ on each
of them, and an ordering $\sigma$ on the union of these sets, that extends
the $\sigma_i$.  This is an example of what is known as the merging of
ordered sets.

Now we will look at $f^vG$.  The main idea is that it can be calculated
from the $f^v T_i$.  Considered in isolation each acyclic shelling
$\sigma_i$ of each $T_i$ contributes $b^{n(i)-1}a$, where $n(i)$ is the
number of vertices on $T_i$.  If we have such a $\sigma_i$ for each
component $T_i$, then the total contribution due to the $\sigma_i$ will be
the sum of all the possible `mergings' of the $b^{n(i)-1}a$.  Here, each
merging of the $\sigma_i$ to $\sigma$ will induce in a natural way a
`merging' of the $b^{n(i)-1}a$.

It should now be clear that $f^vG$ is determined in a straightforward
manner from the $f^vT_i$ of its components.  This follows from a technical
lemma, whose proof is left to the reader.

\begin{lemma}
Suppose $G=G_1 \sqcup G_2$ is a disjoint union of two graphs with only
optional edges.  Then there is a universal bilinear function $\mu$ for
computing $f^vG$.  In other words, the equation
\[
    f^vG = \mu ( f^vG_1 , f^vG_2 )
\]
holds.  Moreover, the coefficients are all nonnegative integers.  (By
induction, there will also be a multilinear $\mu$ for when $G$ is broken
into more than two components $G_i$.)
\end{lemma}

We now come to the main definitions and theorem of this section.

\begin{definition}
Let $T$ be a tree on $n \geq 3$ vertices.  Then the \emph{shelling number}
$s(T)$ of $T$ is the number of ways of shelling $T$ down to a $3$-vertex
tree.  If $T$ has $3$ or fewer vertices, then $s(T)$ is defined to be $1$.
\end{definition}

As a $3$-vertex tree has precisely $4$ acyclic shellings, $s^a(T) = 4
s(T)$ for $n \geq 3$.  For $n=2$, $s^a(T) = 2 s(T)$. For $n=1$,
$s^a(T)=s(T)$.

\begin{definition}
The \emph{shelling number} $s(G)$ of an acyclic graph $G$ is the product
$s(T_1) \dots s(T_r)$ of the shelling numbers of its component trees
$T_i$.  If $G$ contains a cycle, then $s(G)$ is defined to be zero.
\end{definition}

\begin{definition}
Let $G$ be a graph on $n$ vertices.  The \emph{concise flag vector} $f^cG$
is a formal sum of partitions of $n$.  Each subgraph $H$ contributes
$s(H)\pi(H)$ to $f^cG$, where $s(H)$ is as above, and $\pi(H)$ is the
partition of $n$ induced by $H$-connectivity on the vertex set.
\end{definition}

The previous results, taken together, prove the following:

\begin{theorem}
Let $G$ be a graph.  Then the verbose flag vector $f^vG$ is a linear
function of the concise flag vector $f^cG$.  Further, the subgraph and
concise flag vectors are linear functions of each other.
\end{theorem}

In the next section we will show that $f^cG$ is a linear function of
$f^vG$, and thus the concept of a linear function of the flag vector does
not depend on the form of flag vector that is chosen.  We will also show
that the $f^cG$ provide an integral spanning set for the space of formal
sums of partitions.

\section{Basis and dimension}

In this section we will show that the concise flag vector is a linear
function of the verbose flag vector.  This will complete the proof of
their equivalence.  But first we will produce formal sums of graphs, whose
flag vectors provide a basis for formal sums of partitions.

For this we will need some notation.  We will let $[1+1+2]$ denote a
partition of $4$, and $[1+1+2]+[3]$ will denote a formal sum of two
partitions.  If we need to do arithmetic within a partition, we will
enclose it by parentheses, like so $[1 + (n-1)]$.  We will need some
graphs.  We will use $A_n$ and $D_n$ to denote the corresponding
Coxeter-Dynkin diagrams, considered as graphs.  Thus, $A_n$ is $n$
vertices in a row, each joined by an edge to the next.  The graph $D_n$
has a central vertex, from which three arms extend, two of length one, and
one of length $n-3$. We set $D_3$ equal to $A_3$.

Throughout this section $A_n$ and $D_n$ will denote graphs with only
optional edges, considered as formal sums of actual graphs.  Their concise
flag vectors are as follows.

\begin{lemma}
The concise flag vector $f^cA_n$ is given by
\[
    f^c(A_1) = [1] \>, \quad
    f^c(A_2) = [2] \>, \quad
    f^c(A_3) = [3] \>, \quad
    f^c(A_n) = 2^{n-3}[n] \quad \mbox{for $n \geq 3$}
\]
while for $D_n$ we have
\[
    f^c(D_3) = [3] \>, \quad
    f^c(D_4) = 3[4] \>, \quad
    f^c(D_5) = (2\times 2+ 3) [5] = 7 [5]
\]
when $n$ is small, and $f^c(D_n) = (2^{n-2} -1) [n]$ is the general rule.
\end{lemma}

\begin{proof}
The $A_n$ part of this lemma is straightforward.  We simply remove edges
one at a time from either end, until only two edges are left.  Each
additional edge gives another two-fold choice to make.  The $D_n$ part is
not difficult either.  The first edge to be removed can be one of the two
short arms, or can come from the end of the large arm.  Removing a short
arm from $D_n$ turns it into $A_{n-1}$, which we know how to do.  Removing
an edge from the end of the long arm turns $D_n$ into $D_{n-1}$, which we
also know, this time by induction.
\end{proof}

The following result is now obvious.  It can also be seen by comparing the
shellings of $2A_n$ and of $D_n$.

\begin{proposition}
$f^c(2A_n-D_n) = [n]$ for $n \geq 3$.
\end{proposition}

Now suppose that $\pi=[n(1) + \dots + n(r)]$ is a partition of $n$.  It now
follows that the formal expression
\[
    (2A_{n(1)} - D_{n(1)}) \sqcup \dots \sqcup (2A_{n_(r)} - D_{n(r)})
\]
considered as a formal sum of graphs on $n$ vertices has concise flag
vector the given partition $\pi$.  This proves:

\begin{theorem}
The concise flag vectors of graphs on $n$ vertices are an integral
spanning set for the formal sums of partitions of $n$.
\end{theorem}

Recall that $f^vG$ is a linear function of $f^cG$.  We have just shown
that $f^cG$ spans a space of dimension $p(n)$, the number of partitions of
$n$.  If we can show that $f^vG$ also spans a $p(n)$-dimensional space, it
will follow that $f^cG$ is a linear function of $f^vG$.  It is to this
task that we now turn.

We will produce $p(n)$ graphs, or more exactly formal sums of graphs,
whose verbose flag vectors can be seen to be linearly independent. 
Although not providing an integral basis, the simplest graphs for this
purpose are disjoint unions of $A_i$, considered as always in this section
as graphs with only optional edges.

Consider for example $G=A_1\sqcup A_2 \sqcup A_4$.  This is an acyclic
graph with $3$ components and $4$ edges.  From this it follows that for
the $w$-component $f^v_wG$ to be non-zero, $w$ must be a product of $3$
$a$'s and $4$ $b$'s.  However, not for every such product will $f^v_wG$ be
non-zero.  For example, if $w=aabbbba$ then $f^v_wG$ is zero.  This is
because of the leading pair of $a$'s.  It is not possible to begin the
shelling of this $G$ by removing two vertices, without along the way
removing an edge.  (Acyclic shellings only contribute.  If $G$ were $A_1
\sqcup A_1 \sqcup A_5$ it could be done.)

We can now formulate an upper diagonal argument that will prove that the
vectors in our candidate basis are indeed linearly independent.  Let $\pi$
be a partition of $n$, and let $G_\pi$ be the disjoint union of $A_i$
corresponding to $\pi$, thought of as usual as a graph with only optional
edges.  Now let $w(\pi)$ be the first word $w$ in $a$ and $b$, under the
lexicographic order, for which $f^v_w G_\pi$ is nonzero.  The argument of
the previous paragraph shows the following lemma, from which the theorem
follows easily.

\begin{lemma}
If $\pi=[n(1)+ \dots + n(r)]$ is a partition of $n$, and the $n_i$ are in
non-decreasing order, then
\[
    w(\pi) = b^{n(1)-1} a \cdot b^{n(2)-1} a \cdot 
            \ldots \cdot b^{n(r)-1}a
\]
and so the $w(\pi)$ are distinct, as $\pi$ ranges over all partitions of
$n$.
\end{lemma}

\begin{theorem}
The $p(n)$ vectors $f^v G_\pi$, as $\pi$ ranges over the partitions of
$n$, are linearly independent.  Thus, the verbose and the concise flag
vectors are linear functions of each other.
\end{theorem}

\begin{proof}
Let $\pi$ and $\pi'$ be two partitions of $n$.  We have shown that the
matrix $f^v_{w(\pi')} G_\pi$ is upper triangular, when partitions are
ordered lexicographically via $w(\pi)$ and $w(\pi')$.  The diagonal
entries are clearly non-zero, and so the result follows.
\end{proof}

\section{Graphs on $4$ vertices}

It is possible to calculate by hand the flag vectors of all graphs on $4$
vertices.  With a certain amount of care the $3_0$, $3_1$, $3_2$ and $3_3$
notation used for the $3$-vertex graphs can be extended to this case.

A graph on $4$ vertices can have at most $6$ edges.  Let $4_0$ and $4_1$
denote the $4$-vertex graphs with $0$ and $1$ edges respectively.  There
are two $4$-vertex graphs with $2$ edges.  They can be denoted by
$2_1\sqcup2_1$ and $3_2\sqcup1_0$ respectively.  There are three
$4$-vertex graphs with $3$ edges.  One of them is $3_3\sqcup1_0$, while
another is its complement $\overline{3_3\sqcup1_0}$, which is also known
as $D_4$.  The third we will denote by $4_3$, which is also known as
$A_4$.  (Here $A_4$ and $D_4$ are graphs with ordinary, not optional,
edges.)  Finally, the graphs with $4$, $5$ and $6$ edges are the
complements of graphs with $0$, $1$ and $2$ edges.  Thus, we have a
compact notation for every graph on four vertices.

\begin{table}
\begin{center}
\advance\baselineskip 2pt
\begin{tabular}{l|ccccc}
 & $[1+1+1+1]$ & $[2+1+1]$ & $[2+2]$ & $[3+1]$ & $[4]$ \\
\hline
 $4_0$ &                            1&0&0&0&0 \\
 $4_1$ &                            1&1&0&0&0 \\
 $2_1\sqcup2_1$ &                   1&2&1&0&0 \\
 $3_2\sqcup1_0$ &                   1&2&0&1&0 \\
 $4_3=A_4$ &                        1&3&1&2&2 \\
 $3_3\sqcup1_0$ &                   1&3&0&3&0 \\
 $\overline{3_3\sqcup1_0}=D_4$ &    1&3&0&3&3 \\
 $\overline{3_2\sqcup1_0}$ &        1&4&1&5&7 \\
 $\overline{2_1\sqcup2_1}$ &        1&4&2&4&8 \\
 $\overline{4_1}=4_5$ &             1&5&2&8&18 \\
 $\overline{4_0}=4_6$ &             1&6&3&12&36 \\
\end{tabular}
\end{center}
\caption{The concise flag vectors for graphs on $4$ vertices}
\end{table}

\begin{theorem}
The concise flag vectors of the graphs on $4$ vertices are as in Table~1. 
Each such flag vector is a distinct vertex of the convex hull $\Delta(4)$
of these vectors.
\end{theorem}

\begin{proof}
The first statement is left to the reader.  The second is the result of
applying the PORTA convex polytope software package to these values. 
(Alternatively, it could be computed by hand.)
\end{proof}

\section{Flags}

Flags are important in the theory of convex polytopes not only numerically
but also geometrically.  Let $\delta$ be a flag $(\delta_1\subset \dots
\subset \delta_r \subset \Delta)$ on a convex polytope $\Delta$.  From
$\delta$ a flag $\langle\delta\rangle$ of vector spaces can readily be
constructed.  (Each $\langle\delta\rangle_i$ will be the span of the
vectors lying on $\delta_i$.)  Using linear algebra all manner of vector
spaces can be constructed from this flag.  The standard formula
\cite{bib.JD-FL.IHNP,bib.KF.IHTV,bib.RS.GHV} for the middle perversity
intersection homology Betti numbers $h_i\Delta$ of $\Delta$ unwinds
\cite{bib.JF.CPLA} to suggest a construction of a vector space
$\Lambda(\delta,i)$ attached to each flag $\delta$, such that the dimension
of $\Lambda(\delta,i)$ is precisely the contribution $\delta$ makes to
$h_i\Delta$.

That the $h_i\Delta$ are non-negative, at least when $\Delta$ has rational
vertices, is a subtle combinatorial inequality.  When $\delta'$ is
obtained from $\delta$ by removing any one of the $\delta_i$ from
$\delta$, the space $\Lambda(\delta,i)$ is contained in
$\Lambda(\delta',i)$.  This allows the construction of a complex of vector
spaces that is, the author conjectures, exact except at just one
location.  If this is true then the homology at this location will by
construction have dimension $h_i\Delta$.

In much the same way, vector spaces can be attached to each flag on a
graph.  This is an important property, without which the flag vector
concept would in some sense be inadequate. 

\begin{definition}
A \emph{semi-concise flag} $\delta$ of type $w$ on a graph $G$ on $n$
vertices consists of the following.  Write $w$ as $a^{d(1)}b\cdot a^{d(2)}
\ldots$.  To each $d(i)$ associate a $d(i)$-element subset $S_i$ of the
vertices, and to each $b$ associate a vertex $v_i$.  The sets and vertices
are to be disjoint.  Call this a \emph{$w$-partition} of the vertex set
$V$.  In addition, to each $b$ associate an edge from $v_i$ to some vertex
that either lies in $S_j$ or is $v_j$, for some $j>i$.  Such a
configuration is a \emph{semi-concise flag}.
\end{definition}

Now suppose that $\delta$ is a semi-concise flag on $G$, whose word $w$
ends in $ba$.  By permuting the last two vertices another semi-concise
flag $\delta'$ can be obtained.  But $\delta'$ gives no new combinatorial
information, and so $\delta$ and $\delta'$ should be counted together. 
The same argument applies when $\delta$ is a semi-concise flag ending in
$ba^{d(r)}$.  In this case, permuting the ends of this last edge produces
a $\delta'$, that should be counted along with $\delta$.

Now suppose that $\delta$ is a semi-concise flag whose word ends with
$bba$.  This means that we are given two edges that join the last three
vertices.  Now suppose we are given the last three vertices, in no
particular order, and two edges linking them.  For four out of the six
possible permutations of the last three vertices, we will obtain the tail
end of a flag whose word ends in $bba$.  These four flags should be
counted together, because each contains the same information as the
others.

\begin{definition}
A \emph{concise flag} $\delta$, or \emph{flag} for short, is the same as a
semi-concise flag except that one consolidates information as in the
previous two paragraphs.  
\end{definition}

Now let $E$ be the $n$-dimensional vector space that has the vertices
$v_i$ of $G$ as basis vectors.  Each concise flag $\delta$ of $G$ will
induce subspaces in $E$, whose dimensions and mutual dispositions are
determined by the type $w$ of the flag $\delta$.  As in convex polytopes,
one can now seek to interpret linear functions of flag vectors in terms of
linear algebra constructions on the induced subspaces of $E$.  However, at
present we do not have any idea as to which linear functions of the graph
flag vector will play the r\^ole that the $h_i\Delta$ do for polytope flag
vectors.

\section{Summary and conclusions}

We have defined for each graph $G$ on $n$ vertices a flag vector $fG$,
which lies in a space whose dimension is the number $p(n)$ of partitions
of $n$.  This flag vector can be given in verbose, subgraph and concise
forms.  Each is a linear function of the others.  Here is another form for
the flag vector.

\begin{definition}
Let $\Gcal_n$ denote the space of all formal sums of graphs on $n$
vertices, up to isomorphism.  Now define the \emph{nullspace}
$\Zcal=\Zcal(\Gcal_n)$ to be all such formal sums $G$, such that the flag
vector $fG$ is zero.  The \emph{abstract flag vector} $f^aG$ of a graph
$G$ is the residue of $G$, as a vector in $\Gcal_n$, in the quotient space
$\Gcal_n/\Zcal_n$.
\end{definition}

We already know that graphs with an optional cycle will generate elements
of nullspace $\Zcal_n$.  It is not clear whether or not these elements
will span the whole of the $\Zcal_n$.  If not, then finding a pleasant
geometric characterisation of $\Zcal_n$ is an open problem.

Next we turn to $\Delta(n)$, or in other words to the linear inequalities
satisfied by the flag vector.  For $n\leq 4$ each distinct graph gives a
distinct vertex of $\Delta(n)$, or in other words each graph is extremal
for some linear function of the flag vector.  On would like to know if the
same holds for, say, $n \leq 10$.  (This statement is much stronger than
saying the flag vector distinguishes graphs, although of course it is
still considerably weaker than the embedding of graphs in $\Gcal_n$.)

The theory of flag vectors for graphs can be generalised in two ways.  The
first is an extension to hypergraphs.  A graph is a collection $E$ of
$2$-element subsets of the vertex set $V$.  An $i$-graph is a collection
of $i$-element subsets of $V$.  In \cite{bib.JF.QTHGFV}, and more
concisely in \cite{bib.JF.SFV}, a flag vector is defined for such
hypergraphs, and various other objects besides.

The second generalisation is this.  We have been thinking of a graph as a
collection of vertices joined by edges, or in other words, the natural
companions of a graph $G$ on $n$ vertices are the other graphs on $n$
vertices.  We can however think of a graph as a collection of edges joined
at vertices.  When this is done, the number of edges is paramount, and
shelling consists of removing the edges one at a time.

The edge flag vector of $G$ (this paper has hitherto studied the vertex
flag vector) can be defined in the following way.  Each edge has two ends,
and at each end there is a multiplicity.  When in a shelling an edge is
removed, record (a)~the fact that the edge is removed, (b)~the smaller of
the end multiplicities, less one for the edge itself, and (c)~the
difference between the larger and smaller of the end multiplicities.  Use
expressions such as $a + 4 b + 2 c$ to record this data.  (Perhaps
whimsically, $a$ stands for always, $b$ for both, and $c$ for choice.)

Each edge shelling will thus determine a sequence of linear expressions in
$a$, $b$ and $c$.  Form their product, and define the edge flag vector to
be the formal sum of these products, over all edge shellings.  This edge
flag vector has not yet been investigated.

\section*{Acknowledgements}

Thanks are due to Imre Leader and Tomas Slivnik for their interest at
various times in this work, and also to Marge Bayer for doing a PORTA
calculation.

\end{document}